\documentclass[a4paper,11pt]{article}
\usepackage{amssymb,amsmath,amsthm} 
\usepackage{tikz}
\usepackage{comment}

\newtheorem{thm}{\bfseries Theorem}
\newtheorem{lem}[thm]{\bfseries Lemma}

\newtheorem{prop}[thm]{\bfseries Proposition}
\newtheorem{cor}[thm]{\bfseries Corollary}

\newtheorem{prob}[thm]{\bfseries Problem}

\newcommand{\Real}{{\mathbb R}}

\def\CC{{\mathcal C}}

\def\HH{{\mathcal H}}

\def\PP{{\mathcal P}}

\def\SS{{\mathcal S}}

\def\YY{{\mathcal Y}}

\definecolor{leaf}{rgb}{0.2,0.5,0.2}

\newenvironment{prooof}{\medskip
\noindent{\scshape Proof:}}{\quad $\Box$\medskip}

\title{A coding problem for pairs of subsets}

\author{B\'ela Bollob\'as%
\thanks{Department of Pure Mathematics and Mathematical Statistics,
Wilberforce Road, Cambridge CB3 0WB, UK, and
Department of Mathematical Sciences, University of Memphis, Memphis TN 38152, USA.}
\and Zolt\'an F\"uredi%
\thanks{Alfr\'ed R\'enyi Institute of Mathematics, 13--15 Re\'altanoda Street, 1053 Budapest, Hungary. \newline
E-mail: {\tt z-furedi@illinois.edu}.
Research supported in part by the Hungarian National Science Foundation OTKA 104343,
and by the European Research Council Advanced Investigators Grant 267195.}
\and Ida Kantor%
\thanks{Computer Science Institute of Charles University, Malostransk\'e n\'am. 25, 118 00 Praha 1.
Czech Republic.  \newline E-mail: {\tt ida@iuuk.mff.cuni.cz}.
Research supported by GA\v{C}R grant number P201/12/P288 and partially
done while this author visited the R\'enyi Institute.}
\and G. O. H. Katona%
\thanks{R\'enyi Institute, Hungarian Academy of Sciences, Budapest, Re\'altanoda u. 13--15, 1053 Hungary.
 \newline E-mail: {\tt ohkatona@renyi.mta.hu}.
Research was supported by the Hungarian National Foundation OTKA NK104183.
This work was done while this author visited the University of Memphis.}
\and Imre Leader%
\thanks{Department of Pure Mathematics and Mathematical Statistics,
Wilberforce Road, Cambridge CB3 0WB, UK.\enskip E-mail{\tt I.Leader@dpmms.cam.ac.uk}.
This work was done while this author visited the University of Memphis.}
}

\begin{document}
\maketitle

{\small{\em Abstract:}\enskip
Let $X$ be an $n$--element finite set, $0<k\leq n/2$ an integer.
Suppose that $\{ A_1,A_2\} $ and $\{ B_1,B_2\} $ are pairs of disjoint
$k$-element subsets of $X$ (that is, $|A_1|=|A_2|=|B_1|=|B_2|=k$, $A_1\cap
A_2=\emptyset$, $B_1\cap B_2=\emptyset$). Define the distance of these pairs by
$d(\{ A_1,A_2\} ,\{ B_1,B_2\} )=\min \{ |A_1-B_1|+|A_2-B_2|,
|A_1-B_2|+|A_2-B_1|\} $.
This is the minimum number of elements of $A_1\cup A_2$ one has to move to obtain the other pair $\{B_1,B_2\}$.
Let $C(n,k,d)$ be the maximum size of a family of pairs of disjoint $k$-subsets, such
that the distance of any two pairs is at least $d$.}

{\small
Here we establish a conjecture of Brightwell and Katona concerning an asymptotic formula
 for $C(n,k,d)$  for $k,d$ are fixed and $n\to \infty$.
Also, we find the exact value of $C(n,k,d)$ in an infinite number of cases,
 by using special difference sets of integers.
Finally, the questions discussed above are put into a more general context and a
 number of coding theory type problems are proposed.}

\begin{quote}{\small
}
\end{quote}
\vspace{5mm}

\section{The transportation distance}\label{intro}

Let $X$ be a finite set of $n$ elements.
When it is convenient we identify it with the set $[n]:=\{ 1, 2, \dots, n\}$.
The family of the $k$-sets of an underlying set $X$ is denoted by ${X \choose k}$.
For $0<k \leq  n/2$ let $\YY$ be the family of unordered disjoint pairs $\{A_1,A_2\}$
 of $k$-element subsets of $X$ (that is, $|A_1|=|A_2|=k, A_1\cap A_2=\emptyset$).
The {\em transportation distance} or {\em Enomoto-Katona distance} $d$ on $\YY$ is defined by
\begin{equation}
d(\{ A_1,A_2\} ,\{ B_1,B_2\} )=\min \{ |A_1-B_1|+|A_2-B_2|,
 |A_1-B_2|+|A_2-B_1|\} \label{eq:1distdef}.
\end{equation}

In fact, this is an instance of a more general notion.
Whenever $(Z,\rho)$ is a metric space, we can define a metric $\rho^{(s)}$ on $Z^{(s)}$,
 the set of unordered $s$-tuples from $Z$, by
\begin{equation}
\rho^{(s)}(\{x_1,\dots,x_s\},\{y_1,\dots,y_s\})=\min_{\pi\in S_s} \sum_{i=1}^s\rho(x_i,y_{\pi(i)}). \label{assignp}
\end{equation}
It is not hard to verify that $\rho^{(s)}$ satisfies the triangle inequality, i.e., it really is a metric.
The transportation distance defined above is obtained by taking $s=2$,
 $Z$ to be the set of $k$-elements subsets of $X$ and $\rho$ is half of their symmetric difference.

The minimization problem~(\ref{assignp}) (where $\rho$ can be an arbitrary metric) is one of the fundamental combinatorial optimization problems, a so called {\em assignment problem}, a special case of a more general {\em Monge-Kantorovich transportation problem} (see, e.g., the monograph~\cite{villani}).

The transportation distance between finite sets of the same cardinalities is
 one of the interesting  measurements
 among many different ways to define how two sets differ from each other.
In~\cite{tusnady}, Ajtai, Koml\'{o}s and Tusn\'{a}dy considered the assignment problem from a different perspective, and determined with high probability the transportation distance between two sets of points randomly chosen in a unit square.

Since the transportation distance is an important notion, especially
 from the algorithmic point of view,  there are monographs and graduate texts
 about this topic, see, e.g., \cite{villani}.
It is also mentioned in the {\itshape Encyclopedia of Distances}~\cite{Deza}
 as the ``KMMW metric'' (p. 245 in Chapter 14) or as the ``$c$-transportation distance''.
Nevertheless, many combinatorial problems are still unsolved.
The packing of sets in spherical spaces with large transportation distance will be discussed
 in~\cite{fur}.

\section{Packings and codes}\label{S2}

Given a metric space $(Z,\rho)$ and a distance $h > 0$,  the {\em packing number}
 $\delta(Z,\geq h)$
 is the maximum number of elements in $Z$ with pairwise distance at least $h$.

A $(v,k,t)$ packing ${\mathcal P}\subseteq {[v]\choose k}$
  is a family of $k$-sets with pairwise intersections at most $t-1$ (here $v\geq k\geq t\geq 1$).
In other words, every $t$-subset is covered at most once.
Its maximum size is denoted by $P(v,k, t)$.
Obviously,
\begin{equation}\label{eq:3}
P(v,k,t)\leq {v\choose t}/{k \choose t}.
  \end{equation}
If here equality holds then $\mathcal P$ is called a Steiner system
 $S(v,k,t)$, or a {\em $t$-design} of parameters $v,k,t$ and $\lambda =1$
  (for more definitions concerning symmetric combinatorial structures
  esp., difference sets, etc. see, e.g., the monograph by Hall~\cite{Hall}).
More generally, for a set $K$ of integers, a family $\mathcal P$ on $v$ elements
 is called a $(v,K,t)$-design (packing) if every $t$-subset of $[v]$ is contained in exactly
  one (at most one) member of $\mathcal P$ and $|P|\in K$ for every $P\in \mathcal P$.

Determining the packing number is a central problem of Coding Theory,
  it is essentially the same problem as finding the rate of a large-distance
  error-correcting code.

If equality holds in \eqref{eq:3} then every $i$-subset of $[v]$ is
 contained in ${v-i\choose t-i}/{k-i \choose t-i}$ members of $\mathcal P$
  for $i=0,1,\dots, t-1$.
We say that $v, k$, and $t$ satisfy the {\em divisibility conditions}
 if these $t$ fractions are integers.
It was recently proved by Keevash~\cite{Keevash} that for any given $k$ and $t$
 there exists a bound $v_0(k,t)$ such that these trivial
 necessary conditions are also sufficient for the existence of a $t$-design.
\begin{equation}\label{eq:K}
  \text{An } S(v,k,t) \text{ exists if } v,k, \text{ and }t
    \text{ satisfy the divisibility conditions and }v> v_0(k,t).
   \end{equation}
This implies R\"odl's theorem\cite{Rodl}, that for given $k$ and $t$ as $v\to \infty$
\begin{equation}\label{eq:W}
  P(v,k,t)= (1+o(1)) {v\choose t}/{k \choose t}.
   \end{equation}
Even more, \eqref{eq:K} implies that here the error term is only $O(v^{t-1})$.
The case $t=2$ was proved much earlier by Wilson~\cite{Wilson}.
For this case he also proved the following more general version.
For a finite $K$ there exists a bound $v_0(K,2)$ such that for $v> v_0(K,2)$
\begin{equation}\label{eq:W2}
  \text{a } (v,K,2) \text { design exists if }  v \text{ and }K
    \text{ satisfy the generalized divisibility conditions},
   \end{equation}
namely, g.c.d.$( {k\choose 2}: k\in K) $ divides ${v\choose 2}$
 and g.c.d.$( k-1: k\in K )$ divides $v-1$.

\section{Packing pairs of subsets}\label{S3}
In this paper, we concentrate on the space $\YY$ of pairs of {\em disjoint} $k$-subsets.
We say that a set ${\cal C}\subset \YY$ of such pairs is a $2$-$(n,k,d)$--{\em code} if
the distance of any two elements is at least $d$.
Let $C(n,k,d)$ be the maximum size of a $2$-$(n,k,d)$-code.
Enomoto and Katona in~\cite{Eno} proposed the problem of determining $C(n,k,d)$.
For the origin of the problem see \cite{DKS}.
Connections to Hamilton cycles in the Kneser graph $K(n,k)$ are discussed in~\cite{Kat}.
The problem makes sense only when $d\leq 2k\leq n$.
It is obvious, that a maximal $2$-$(n,k,1)$ code consists of all the pairs,
 $C(n,k,1)=|\YY|= {1\over 2}{n\choose k}{n-k\choose k }$.
A $2$-$(n,k,2k)$ code consists of mutually disjoint $k$-sets, hance
  $C(n,k,2k)=\lfloor n/2k\rfloor$.

In Section~\ref{sec2} we present a method for the determination the exact value of $C(n, k, 2k-1)$
 for infinitely many $n$.
However, we were able to complete the cases $k=2,3$ only, the cases of pairs and triple systems.
\begin{thm}\label{th:23}
If $n\equiv 1\bmod 8$ and $n> n_0$  then $C(n,2,3)={n(n-1)\over 8}$.\\
If $n\equiv 1,19\bmod 342$ and $n> n_0$ then $C(n,3,5)={n(n-1)\over 18}$.
\end{thm}

The following theorem was proved in~\cite{BK}.
Let $d\leq 2k\leq n$ be integers. Then
\begin{equation}\label{eq:upp}
 C(n,k,d)
  \leq {1\over 2}{n(n-1)\cdots (n-2k+d)\over k(k-1)\cdots
\lceil {d+1\over 2}\rceil \cdot k(k-1)\cdots \lfloor {d+1\over 2}\rfloor }.
  \end{equation}
Quisdorff~\cite{Quis} gave a new proof and using ideas from classical
 coding theory he significantly improved the upper bound for small values of $n$
 (for $n\leq 4k$).
For completeness, in Section~\ref{sec4} we reprove \eqref{eq:upp} in an even more streamlined way.

Concerning larger values of $n$ one can build a $2$-$(n,k,d)$ code
 from smaller ones using the following observation.
If $|(A_1\cup A_2) \cap (B_1\cup B_2)|\leq 2k-d$ holds for the disjoint pairs
 $\{ A_1,A_2\}\in \YY$,  $\{ B_1,B_2\}\in \YY$ then
 $d(\{ A_1,A_2\}, \{ B_1,B_2\})\geq d$.
Take a $(2k-d+1)$-packing $\PP$ on $n$ elements and choose a $2$-$(|P|,k,d)$-code on
 each members $P\in \PP$. We obtain
\begin{equation}\label{eq:genP}
 \sum_{P\in \PP} C(|P|,k,d) \leq  C(n,k,d).
  \end{equation}
This gives
\begin{equation}\label{eq:genP9}
 P(n,p,2k-d+1) C(p,k,d) \leq  C(n,k,d).
  \end{equation}
Fix $p$ (and $k$, $t$ and $d$) then R\"odl's theorem \eqref{eq:W} gives
$(1+o(1)){n \choose 2k-d+1}{p\choose 2k-d+1}^{-1}C(p,k,d) \leq  C(n,k,d)$.
Rearranging we get, that the sequence
  $C(n,k,d)/{n \choose 2k-d+1}$ is essentially nondecreasing in $n$,
  for any fixed  $p$ (and $k$, $t$ and $d$)
\begin{equation*}\label{eq:genP10}
  C(p,k,d)/{p\choose 2k-d+1} \leq (1+o(1)) C(n,k,d)/{n \choose 2k-d+1}.
  \end{equation*}
Since, obviously, $C(2k,k,d)\geq 1$ we obtain that
$\lim_{n \to \infty} C(n,k,d)/{n \choose 2k-d+1}$
exists, it is positive, it equals to its supremum, and finite by~\eqref{eq:upp}.

It was conjectured (\cite{BK}, Conjecture 8) that the upper estimate~(\ref{eq:upp}) is asymptotically sharp.
We prove this conjecture in Section~\ref{sec3}.

\begin{thm}\label{th:2}
$$\lim _{n\rightarrow \infty} {C(n,k,d)\over n^{2k-d+1}}= {1\over 2}{1\over k(k-1)\cdots
\lceil {d+1\over 2}\rceil \cdot k(k-1)\cdots \lfloor {d+1\over 2}\rfloor }.$$
\end{thm}

\section{The case $d=2$, the exact values of $C(n, k, 2)$}\label{sec:d=2}

Besides the cases mentioned in the previous Section
 (the cases $d=1$, $d=2k$ and $(k,d)\in \{ (2,3),(3,5)\}$)
 we can solve one more case easily, namely if $d=2$.
Since $C(2k,k,2)]=|\YY|=\frac{1}{2}{2k\choose k}$ the construction~\eqref{eq:genP9} gives
 $P(n,2k,2k-1)\frac{1}{2}{2k\choose k}\leq C(n,k,2)$.
Then the recent result of Keevash~\eqref{eq:K} gives the lower bound in the following
 Proposition.
The upper bound follows from~\eqref{eq:upp}.
 \begin{prop}\label{prop4} $C(n,k,2)={n\choose 2k-1}\frac{1}{4k}{2k\choose k}$ for all
  $n> n_0(k)$ whenever the divisibility conditions of \eqref{eq:K} hold. \qed
   \end{prop}

\section{The case $d=2k-1$, the exact values of $C(n, k, 2k-1)$}\label{sec2}
The distance $\delta (a,b)$ of two integers $\bmod{\ m}$ $(1\leq a, b\leq m)$ is
defined by
$$\delta (a,b)= \min \{ |b-a|, |b-a+m|\} .$$
(Imagine that the integers $1,2,\ldots ,m$ are listed
around the cirle clockwise uniformly. Then $\delta (a,b)$ is the smaller distance around the
circle from $a$ to $b$.) $\delta (a,b)\leq {m\over 2}$ is trivial.
Observe that $b-a\equiv d-c \bmod m$ implies $\delta (a,b)=\delta (c,d)$.

We say that the pair $S=\{ s_1,\ldots ,s_k\}$,
 $T=\{ t_1,\ldots ,t_k\} \subset \{ 1,\ldots ,m\} $ of disjoint sets is {\it antagonistic}
$\bmod{\ m}$ if

(i)\,\, all the $k(k-1)$ integers $\delta (s_i,s_j)$ $(i\not= j)$ and
$\delta (t_i,t_j)$  $(i\not= j)$ are different,

(ii)\, the $k^2$ integers $\delta (s_i,t_j)$  $(1\leq i,j\leq k)$ are all different and

  (iii) $\delta (s_i,t_j)\not= {m\over 2}$  $(1\leq i,j\leq k).$\newline
If there is a pair of disjoint antagonistic $k$-element subsets $\bmod{\ m}$ then $2k^2+1\leq m$
must hold by (ii) and (iii).

\begin{prob}\label{prob:4}
Is there a pair of disjoint, antagonistic $k$-element sets $\bmod{\ 2k^2+1}$?
\end{prob}

We have an affirmative answer only in three cases.

\begin{prop}\label{prop:5}
There is a pair of disjoint, antagonistic $k$-element sets $\bmod{\ 2k^2+1}$ when $k=1,2,3.$
\end{prop}

\begin{prooof} We simply give such $k$-element sets in these cases.
It is easy to check that they satisfy the conditions.

$k=1$: $S=\{ 1\} , T=\{ 2\} .$

$k=2$: $S=\{ 1,8\} , T=\{ 2,3\} .$

$k=3$: $S=\{ 1,5,19\} , T=\{ 2,13,15\} .$
\end{prooof}

\begin{lem}\label{le:6}
If there is a pair of disjoint, antagonistic $k$-element sets $\bmod{\ \ m}$ then
\newline
$C(m,k,2k-1)\geq m$.
\end{lem}

\begin{prooof} Let $(S,T)$ be the antagonistic pair. The shifts
$S(u)=\{ a+u \bmod m: s\in S\} , T(u)=\{ s+u \bmod m: s\in T\}(0\leq u<m)$ will
serve as pairs of disjoint subsets of $X$.

Suppose that $S(u)$ and $S(v)\ (u\not= v)$ have two elements in common:
$s_1+u=s_2+v\not= s_3+u=s_4+v$ where $s_1, s_2, s_3, s_4\in S, (s_1,s_2)\not= (s_3,s_4) $.
The difference is $s_1-s_2=s_3-s_4$ contradicting (i). One can prove in the same
way that $T(u)$ and $T(v)\ (u\not= v)$ and $S(u)$ and $T(v)$, respectively, have
at most one element in common. In other words the intersection of any pair from
the sets $S(u),T(u),S(v),T(v)$ has at most one element.

Suppose now that both $S(u)\cap S(v)$ and $T(u)\cap T(v)$ are non-empty for some
$u\not= v$. Then $s_1+u=s_2+v, t_1+u=t_2+v$ holds for some $s_1, s_2\in S,
t_1, t_2\in T$. This leads to $v-u=s_1-s_2=t_1-t_2$, contradicting (i), again.

Finally, suppose that both $S(u)\cap T(v)$ and $T(u)\cap S(v)$ are non-empty for some
$u\not= v$. Then $s_1+u=t_1+v, t_2+u=s_2+v$ is true  for some $s_1, s_2\in S,
t_1, t_2\in T$. Here $v-u=s_1-t_1=t_2-s_2$ is obtained, contradicting either (ii)
or (iii) (the latter one, if $s_1-t_1=t_1-s_1$ is obtained).

This proves that the distance of the pairs $(S(u), T(u))$ and $(S(v), T(v))\ (u\not= v)$
is at least $2k-1$.
\end{prooof}

\begin{cor}\label{cor:7}
Suppose that there is Steiner family  ${\cal S}(n,2k^2+1,2)$ and
a disjoint, antagonistic pair of $k$-element subsets $\bmod{\ 2k^2+1}$ then
$$C(n,k,2k-1)={n(n-1)\over 2k^2}.$$
\end{cor}

\begin{prooof}
The upper bound $C(n,k,2k-1)\leq n(n-1)/ 2k^2$ is a corollary of~(\ref{eq:upp}).

The lower estimate is obtained from \eqref{eq:genP9}. By Lemma~\ref{le:6} one can choose
$2k^2+1$ pairs of disjoint $k$-subsets with distance $2k-1$ in a set of $2k^2+1$ elements.
This can be done in each of the members of ${\cal S}(n,2k^2+1,2)$. Since the members have at
most one common element, the distance of two pairs in distinct members of
${\cal S}(n,2k^2+1,2)$ will have distance at least $2k-1$. Therefore all the
$$|{\cal S}(n,2k^2+1,2)| (2k^2+1)={{n\choose 2}\over {2k^2+1\choose 2}}(2k^2+1)={n(n-1)\over 2k^2}$$
pairs have distance at least 1.
\end{prooof}

\noindent
{\sc Proof} of Theorem~\ref{th:23}. \quad
We only need lower bounds, i.e., constructions.
The case $k=3$ follows from Wilson's theorem~\eqref{eq:K} of the existence of $S(n,19,2)$,
 Proposition~\ref{prop:5} and Corollary~\ref{cor:7}.

Similarly, the case $k=2$ for $n\equiv 1,9\bmod 72$ follows in the same way using
 Steiner systems $S(n,9,2)$ and the fact $C(9,2,3)=9$ from Corollary~\ref{cor:7}.
However, one can see that $C(17,3,2)=34$ and then the results follows from
 Wilson's theorem~\eqref{eq:W2} of the existence of $S(n,\{9,17 \},2)$ for all
 large $n\equiv 1 \bmod 8$ and construction~\eqref{eq:genP}.

The construction for $C(17,2,3)$ is similar to the proof of Lemma~\ref{le:6}.
The 9 pairs there are defined as $\{ \{x+1,x+8\}, \{x+2 ,x+3\}\}: x\in Z_9\}$.
These correspond to a perfect edge decomposition of $K_9$ into $C_4$'s with side lengths
 $1,3,4$, and $2$.
For $n=17$ we take the pairs
  $\{ \{x,x+7\}, \{x+2,x+6\}\}: x\in Z_{17}\}$ and $\{ \{y,y+11\}, \{y+7 ,y+8\}\}: y\in Z_{17}\}$
  which correspond to $C_4$'s of side lengths $(2,5,1,6)$ and $(7,4,3,8)$, respectively.
\qed

Note that the method gives that $C(n,1,1)={n(n-1)\over 2}$ when $n \equiv 1,3 \bmod 6$.
This, however, is trivial for all $n$.

\section{A new proof of the upper estimate}\label{sec4}
The upper estimate in~(\ref{eq:upp}) was proved in~\cite{BK}.
We give a new, more illuminating proof here.

Given a pair $\{A,B\}$ of disjoint $k$-element sets let ${\cal P}(\{A,B\},u,v)$ denote the
family of pairs $\{U,V\}$ where $|U|=u, |V|=v$ and $U\subseteq A, V\subseteq B$ or vice versa.
We have
$$ |{\cal P}(\{A,B\},u,v)|=2{k\choose u}{k\choose v}.$$
Suppose first $u<v$.
Then the total number of pairs $\{U,V\}, U\cap V=\emptyset , |U|=u, |V|=v$ in an $n$-element set is
$${n\choose u}{n-u\choose v}.$$

Let $\{A_1,B_1\},\{A_2,B_2\}$ be two pairs with distance at least $d$, and $u < v$ be two nonnegative integers such that $u+v=2k-d+1$. By definition~\eqref{eq:1distdef},
$ {\cal P}(\{A_1,B_1\},u,v)$ and $ {\cal P}(\{A_2,B_2\},u,v)$ are disjoint.
We have
\begin{equation}
C(n,k,d)\leq {{n\choose u}{n-u\choose v}\over 2{k\choose u}{k\choose v}}\label{eq25}
=\frac{n(n-1)\dots (n-2k+d)}{2k(k-1)\dots (k-u+1)k(k-1)\dots (k-v+1)}
\end{equation}
for every pair $u,v$ that satisfies the above requirements.
If $u=v$, then equality~(\ref{eq25}) holds by similar arguments.

The numerator does not depend on $u$, and the denominator is maximized when $u$ and $v$ are as close as possible, i.e., for
$u=2k-\lceil \frac{d-1}{2}\rceil$ and $v=2k-\lfloor \frac{d-1}{2}\rfloor$.
Substituting these values, we obtain the upper estimate in~(\ref{eq:upp}). \qed

\section{Nearly perfect selection}\label{sec3}

Let $\mathcal{W}$ be the family of pairs $\{U,V\}$ such that $U,V\subseteq [n]$, $U\cap V=\emptyset$,
 and $|U|+|V|=2k-d+1$ holds.
Note that $|\mathcal{W}|=\frac{1}{2} \sum_{0\leq u\leq 2k-d+1}{n \choose u}{n-u \choose (2k-d+1)-u}$.
For a pair $\{A,B\}$ of disjoint $k$-element sets, let ${\cal P}(\{A,B\})$ denote the
family of pairs $\{U,V\}\in \mathcal{W}$ for which $U\subseteq A$ and $V\subseteq B$, or vice versa.

\begin{lem}\label{dist}
 $d(\{ A_1,B_1\} ,\{ A_2,B_2\}) \leq d-1$ holds if and only if
${\cal P}(\{A_1,B_1\})\cap {\cal P}(\{A_2,B_2\})\neq \emptyset.$
\end{lem}

\begin{prooof} Suppose that $\{U,V\}\in {\cal P}(\{A_1,B_1\})\cap {\cal P}(\{A_2,B_2\})$,
 say $U\subset A_1\cap A_2$ and $V\subset B_1\cap B_2$.
Then $|A_1-A_2|\leq k-|U|,
|B_1-B_2|\leq k-|V|$ imply $|A_1-A_2|+|B_1-B_2|\leq 2k-|U|-|V|=d-1$ proving the statement. The other case is analogous.

Conversely, if the distance is at most $d-1$ then either $|A_1-A_2|+|B_1-B_2|\leq d-1$ or $|A_1-B_2|+|B_1-A_2|\leq d-1$
must hold. Suppose that the first one is true. Then $|A_1\cap A_2|+|B_1\cap B_2|\geq 2k-d+1$ follows. Take $U=A_1\cap A_2$
and a $V\subseteq B_1\cap B_2$ such that $|V|=2k-d+1-|U|$. Then ${\cal P}(\{A_1,B_1\})\cap {\cal P}(\{A_2,B_2\})\not=
\emptyset $ holds, as claimed.
\end{prooof}

We can view the sets ${\cal P}(\{A,B\})$ as the edges of a hypergraph on the vertex set $\mathcal{W}$.
Let us call this hypergraph $\cal{H}$.
Then a $2$-$(n,k,d)$-code corresponds to a {\em matching} in $\cal{H}$.

In his celebrated paper~\cite{Rodl}, R\"{o}dl established~\eqref{eq:W} in the following way.
He viewed the $t$-element sets as vertices of a $\binom{k}{t}$-uniform hypergraph $\HH_n$
 whose edges correspond to the $k$-element subsets of $[n]$.
Equality~\eqref{eq:W} is in fact a statement about the existence of an almost perfect matching in $\HH_n$.
Using the same key proof idea, a powerful generalization by Frankl and R\"{o}dl~\cite{FR}
guarantees the existence of almost perfect matchings in hypergraphs satisfying certain more general conditions.
Various generalizations and stronger versions versions were later proved, e.g., by Pippenger and Spencer~\cite{PS}.

A function $t:E(\HH)\rightarrow \Real$ is a {\em fractional matching} of the hypergraph $\HH$ if
$\sum_{e\in E(\HH); x\in e}t(e)\leq 1$ holds for every vertex $x\in V(\HH)$.
The {\em fractional matching number}, denoted $\nu^*(\HH)$ is the maximum of $\sum_{e\in E(\HH)}t(e)$ over all fractional matchings. If $\nu(\HH)$ denotes the maximum size of a matching in $\HH$, then clearly
\[
 \nu(\HH)\leq \nu^*(\HH).
\]
Kahn~\cite{Kahn} proved
that under certain conditions, asymptotic equality holds.
Both the hypotheses and the conclusion are in the spirit of the Frankl--R\"{o}dl theorem.

Given a hypergraph $\HH$ with vertex set $[n]$, a fractional matching $t$ and a subset $W\subseteq [n]$,
 define $\bar{t} (W)=\sum_{W\subseteq e\in E(\HH)} t(e)$ and $\alpha(t)=\max\{ \bar{t}(\{x,y\}):x,y\in V(\HH), x\neq y\}.$
In other words, $\alpha(t)$ is a fractional generalization of the codegree.
Let $t(\HH)$ denote $\sum_{e\in E(\HH)}t(e)$.
We say that $\HH$ is $s$-{\em bounded} if each of its edges has size at most $s$.

\begin{thm}[\cite{Kahn}]\label{kahnthm}
For every $s$ and every $\varepsilon>0$ there is a $\delta$ such that whenever $\cal{H}$ is an $s$-bounded hypergraph and $t$ a fractional matching with $\alpha(t)<\delta$, then
\[
 \nu(\HH)>(1-\varepsilon) t(\HH).
\]
\end{thm}

\noindent{\sc Proof} of Theorem~\ref{th:2}. \quad
In the light of Lemma~\ref{dist} it suffices to verify the conditions of Theorem~\ref{kahnthm}
 and to produce a fractional matching $t$ of the hypergraph $\HH$ of the desired size.

Define a constant weight function $t: E(\HH)\rightarrow \mathbb{R}$ by
\[ t(e)=\frac{\lceil\frac{d-1}{2}\rceil! \lfloor\frac{d-1}{2}\rfloor!}{n^{d-1}}.\]
For a vertex $x=\{U,V\}\in \mathcal W$ with $|U|=u$ and $|V|=v$ we have
\[ \deg(\{U,V\})=\binom{n-u-v}{k-u} \binom{n-k-v}{k-v}\leq \frac{n^{d-1}}{(k-u)!(k-v)!}\leq \frac{n^{d-1}}{\lceil\frac{d-1}{2}\rceil! \lfloor\frac{d-1}{2}\rfloor!}\]
hence $t$ is indeed a fractional matching. 
Note that $t(\HH)$ is is asymptotically equal to the quantity in the statement of the Theorem~\ref{th:2}.

The hypergraph $\HH$ is not regular
 but $s$-bounded with $s=\frac{1}{2}\sum_u {k \choose u}{k \choose (2k-d+1)-u}$.
Here $s$ does not depend on $n$.
For $x,y\in V(\HH)=\mathcal{W}$ let $\deg(x,y)$ denote the codegree of $x=\{U,V\}$ and $y=\{U',V'\}$,
i.e., the number of hyperedges ${\cal P}(\{A,B\})$ that contain both $x$ and $y$.
If $U\cup V=U'\cup V'$ (they partition the same $(2k-d+1)$-element set) then the codegre
 $\deg(x,y)=0$.
Otherwise, $|U\cup U'\cup V\cup V'|\geq 2k-d+2$ and $(U\cup U'\cup V\cup V')\subset (A\cup B)$ imply
 that
\[
 \deg (\{U,V\},\{U',V\})=O(n^{d-2}).
\]
Hence $
 \alpha(t)=\deg(\{U,V\},\{U',V\})\cdot t(e)=o(1)$
 and Kahn's theorem completes the proof.  \qed

\section{$s$-tuples of sets, $q$-ary codes}
Let $\YY^{(s)}$ be the family of $s$-tuples of pairwise disjoint $k$-element subsets of $[n]$. A natural definition of a metric on~$\YY^{(s)}$ was already mentioned in the introduction, in equation~(\ref{assignp}). With $\rho$ being half the symmetric difference, the distance is defined as
\[
 \rho^{(s)}(\{A_1,\dots,A_s\},\{B_1,\dots,B_s\})=
\min_{\pi\in S_s} \sum_{i=1}^s |A_i\setminus B_{\pi(i)}|.
\]
Let $C_s(n,k,d)$ denote the maximum size of a subfamily $\SS$ of $\YY^{(s)}$ such that any two elements in~$\SS$ have distance at least $d$. The proofs presented in Sections~\ref{sec3} and~\ref{sec4} can be easily adapted to determining $C_s(n,k,d)$, as well. The proof of the lower and the upper bounds in Theorem~\ref{tuples} is completely analogous to the proofs of inequality~(\ref{eq:upp}) and Theorem~\ref{th:2}.

\begin{thm}\label{tuples}
$$\lim _{n\rightarrow \infty} {C_s(n,k,d)\over n^{sk-d+1}}= {1\over s!}
\frac{{\lceil \frac{d-1}{s}\rceil!} {\lceil \frac{d-2}{s}\rceil!}\dots {\lceil \frac{d-s}{s}\rceil!}}{(k!)^s}. \quad  \qed
  $$
\end{thm}

Let $\YY_q$ be the set of $q$-ary vectors of length $n$ and weight $k$ (weight is the number of nonzero entries). Let $A_q(n,d,k)$ be the maximum size of a subset $\CC\subseteq \YY_q$ such that $\rho'(u,v)\geq d$ whenever $u,v\in \CC$. Here $\rho'$ is the Hamming distance.

With a slightly more technical proof along the same lines, the following can be proven.


\begin{thm}
 Fix $q\geq 2$, $k$ and $d$. If $d$ is odd, then, as $n\rightarrow \infty$,
\[
 A_q(n,d,k)\sim \frac{n^{k-\frac{d-1}{2}} (q-1)^{k-\frac{d-1}{2}} \left(\frac{d-1}{2}\right)!}{k!}.
\]
If $d\geq 2$ is even, then, as $n\rightarrow \infty$,
\[
 A_q(n,d,k)\sim \frac{n^{k-\frac{d}{2}+1} (q-1)^{k-\frac{d}{2}+1} \left(\frac{d}{2}-1\right)!}{k!}.\quad \qed
\]
\end{thm}

To use random methods constructing codes is not a new idea.
The best known general bounds for the {\em covering radius} problems are obtained in this way, see,
e.g.,~\cite{fur203,KSV}.

We can also consider pairs (or more generally $s$-tuples) of $q$-ary vectors of weight $k$.
For simplicity, we will only state the results for pairs here. 
Define the set $\YY^{(2)}_q$ of pairs $\{u,v\}$ of vectors such that
\begin{itemize}
 \item $u,v\in \{0,1,\dots,q-1\}^n$
\item each of $u$ and $v$ has exactly $k$ nonzero entries
\item the supports of $u$ and $v$ are disjoint (i.e. $u_i=0$ for all $i$ such that $v_i\neq 0$, and $v_i=0$ for all $i$ such that $u_i\neq 0$).
\end{itemize}

Define the distance between these pairs by
\[\delta(\{u,v\},\{w,z\})=\min\{\rho'(u,w)+\rho'(v,z),\rho'(u,z)+\rho'(v,w)\}\]
where $\rho'$ is again the Hamming distance.

In the following, $A_q^2(n,d,k)$
will denote the maximum size of a subset $\CC\subseteq \YY^{(2)}_q$ such that $\delta(\{u,v\},\{w,z\})\geq d$ for any pair $\{u,v\},\{w,z\}$ of members of $\CC$.

\begin{thm}\label{2}
 Fix $q$, $d$ and $k$. If $d$ is odd and $q\geq 3$, then, as $n\rightarrow \infty$,
\[A_q^2(n,d,k)\sim \frac{1}{2}\cdot \frac{n^{2k-\frac{d-1}{2}}\cdot (q-1)^{2k-\frac{d-1}{2}} \cdot \lfloor\frac{d-1}{4} \rfloor! \lceil\frac{d-1}{4}\rceil!}{(k!)^2}.\]
If $d\geq 2$ is even and $q\geq 2$, then, as $n\rightarrow \infty$,
\[A_q^2(n,d,k)\sim \frac{1}{2}\cdot \frac{n^{2k-\frac{d}{2}+1}\cdot (q-1)^{2k-\frac{d}{2}} \cdot \lfloor\frac{d}{4} \rfloor! \left(\lceil\frac{d}{4}\rceil -1\right)!}{(k!)^2}. \quad \qed\]
\end{thm}

The distance $\delta$ used here is twice the distance defined in Section~\ref{intro}, hence the apparent inconsistency of this result for $q=2$ with Theorem~\ref{th:2}.

For $q=2$ and $d$ odd we have $A_q(n,d,k)=A_q(n,d+1,k)$.

\section{Open problems}

We believe that for an arbitrary pair of $k$ and $d$, there
are infinitely many $n$'s with equality in inequality~\eqref{eq:upp}.

\section{Further developments}

Let us note that since announcing the first version of the present paper
  Theorem~\ref{th:23} has been greatly extended by Chee, Kiah, Zhang and Zhang~\cite{CKZZ}.
They determined the exact value of $C(n, 2, d)$ completely, and for any fixed $k$
  the exact value of $C(n, k, 2k-1)$ for all $n>n_0(k)$
 satisfying either $n = 0 \mod k$ or
 $n = 1 \mod k$ and $n(n-1) = 0 \mod 2k^2$.
Their proofs are different: they use more design theory.
However, our Section~\ref{sec2} is still interesting for its own sake and Problem~\ref{prob:4} is still open.

\medskip\noindent
{\em Acknowledgements.} \enskip The authors are very grateful for the helpful remarks of the referees.


\begin{thebibliography}{10}

\bibitem{tusnady}
{\sc M.~Ajtai, J.~Koml{\'o}s, and G.~Tusn{\'a}dy}, {\em On optimal matchings},
  Combinatorica, 4 (1984), pp.~259--264.

\bibitem{BK}
{\sc G.~Brightwell and G.~O.~H. Katona}, {\em A new type of coding problem},
  Studia Sci. Math. Hungar., 38 (2001), pp.~139--147.

\bibitem{CKZZ}
{\sc Yeow Meng Chee, Han Mao Kiah, Hui Zhang, and Xiande Zhang},
{\em Optimal codes in the Enomoto-Katona space},
  Combinatorics, Probability and Computing, to appear.
  (Preliminary version in Proc. IEEE Intl. Symp. Inform. Theory. IEEE, 2013.)

\bibitem{DKS}
{\sc J.~Demetrovics, G.~O.~H. Katona, and A.~Sali}, {\em Design type problems
  motivated by database theory}, J. Statist. Plann. Inference, 72 (1998),
  pp.~149--164.
\newblock R. C. Bose Memorial Conference (Fort Collins, CO, 1995).

\bibitem{Deza}
{\sc M. M. Deza and E. Deza},
{\em Encyclopedia of Distances},
\newblock Springer, 2nd ed. 2013. 

\bibitem{Eno}
{\sc H.~Enomoto and G.~O.~H. Katona}, {\em Pairs of disjoint {$q$}-element
  subsets far from each other}, Electron. J. Combin., 8 (2001), Research
  Paper 7, 7 pp. (electronic).
\newblock In honor of Aviezri Fraenkel on the occasion of his 70th birthday.

\bibitem{FR}
{\sc P.~Frankl and V.~R{\"o}dl}, {\em Near perfect coverings in graphs and
  hypergraphs}, European J. Combin., 6 (1985), pp.~317--326.

\bibitem{fur}
{\sc Z. F\"uredi},
  {\em Packings of sets in spherical spaces with large transportation distance}, in preparation.

\bibitem{fur203}
{\sc Z. F\"uredi and J-H. Kang},
Covering the $n$-space by convex bodies and its chromatic number,
{\it Discrete Mathematics} {\bf 308} (2008), 4495--4500.

\bibitem{Hall}
{\sc M. Hall}, {\em Combinatorial Theory, Second Edition}, Wiley-Interscience, 1998.

\bibitem{Kahn}
{\sc J.~Kahn}, {\em A linear programming perspective on the
  {F}rankl-{R}\"odl-{P}ippenger theorem}, Random Structures Algorithms, 8
  (1996), pp.~149--157.

\bibitem{Kat}
{\sc G.~O.~H. Katona}, {\em Constructions via {H}amiltonian theorems}, Discrete
  Math., 303 (2005), pp.~87--103.

\bibitem{Keevash}
{\sc P. Keevash},
{\em The existence of designs}, arxiv.org 1401.3665.

\bibitem{KSV}
{\sc M. Krivelevich, B. Sudakov, and Van H. Vu},
{\em Covering codes with improved density},
IEEE Trans. Inform. Theory 49 (2003), no. 7, 1812--1815.


\bibitem{PS}
{\sc N.~Pippenger and J.~Spencer}, {\em Asymptotic behavior of the chromatic
  index for hypergraphs}, J. Combin. Theory Ser. A, 51 (1989), pp.~24--42.

\bibitem{Quis}
{\sc J\"orn Quistorff}, {\em New upper bounds on Enomoto--Katona's coding type problem},
Studia Sci. Math. Hungar. 42 (2005), pp.~61--72.


\bibitem{Rodl}
{\sc V.~R{\"o}dl}, {\em On a packing and covering problem}, European J.
  Combin., 6 (1985), pp.~69--78.

\bibitem{villani}
{\sc C.~Villani}, {\em Topics in optimal transportation}, vol.~58 of Graduate
  Studies in Mathematics, American Mathematical Society, Providence, RI, 2003.

\bibitem{Wilson}
{\sc R.~M. Wilson}, {\em An existence theory for pairwise balanced designs.
  {II}. {T}he structure of {PBD}-closed sets and the existence conjectures}, J.
  Combinatorial Theory Ser. A, 13 (1972), pp.~246--273.

\end{thebibliography}
\end{document}